\documentclass[10pt]{amsart}

\usepackage[latin1]{inputenc}
\usepackage{amsmath}
\usepackage{amsfonts}
\usepackage{amssymb}
\usepackage[mathscr]{eucal}
\usepackage{diagrams}
\usepackage{xy,color,graphicx}
\usepackage{hyperref}
\xyoption{all}

\newcommand{\wis}[1]{{\text{\em \usefont{OT1}{cmtt}{m}{n} #1}}}

\newcommand{\C}{\mathbb{C}}
\newcommand{\Z}{\mathbb{Z}}
\newcommand{\vtx}[1]{*+[o][F-]{\scriptscriptstyle #1}}

\newtheorem{theorem}{Theorem}

\title{High-dimensional representations of the 3-component loop braid group}
\author{Lieven Le Bruyn} 
\address{Department of Mathematics, University of Antwerp \\ 
 Middelheimlaan 1, B-2020 Antwerp (Belgium) \\ {\tt lieven.lebruyn@uantwerpen.be}}

\begin{document}
\sloppy

\maketitle

\begin{abstract}
In a recent paper \cite{Bruillard} it is shown that irreducible representations of the three string braid group $B_3$ of dimension $\leq 5$ extend to representations of the three component loop braid group $LB_3$. Further, an explicit $6$-dimensional irreducible $B_3$-representation is given not allowing such an extension. 

In this note we give a necessary and sufficient condition, in all dimensions, on the components of irreducible representations of the modular group $\Gamma$ such that sufficiently general representations extend to $\Gamma \ast_{C_3} S_3$. As a consequence, the corresponding irreducible $B_3$-representations do extend to $LB_3$.
\end{abstract}

\section{The strategy}

The $3$-component loop braid group $LB_3$ encodes motions of $3$ oriented circles in $\mathbb{R}^3$. The generator $\sigma_i$ ($i=1,2$) is interpreted as passing the $i$-th circle under and through the $i+1$-th circle ending with the two circles' positions interchanged. The generator $s_i$ ($i=1,2$) simply interchanges the circles $i$ and $i+1$. For physical motivation and graphics we refer to the paper by John Baez, Derek Wise and Alissa Crans \cite{Baez}. The defining relations of $LB_3$ are:
\begin{enumerate}
\item{$\sigma_1 \sigma_2 \sigma_1 = \sigma_2 \sigma_1 \sigma_2$}
\item{$s_1 s_2 s_1 = s_2 s_1 s_2$}
\item{$s_1^2 = s_2^2 = 1$}
\item{$s_1 s_2 \sigma_1 = \sigma_2 s_1 s_2$}
\item{$\sigma_1 \sigma_2 s_1 = s_2 \sigma_1 \sigma_2$}
\end{enumerate}
Note that $(1)$ is the defining relation for the $3$-string braid group $B_3$, $(2)$ and $(3)$ define the symmetric group $S_3$, therefore the first three relations describe the free group product $B_3 \ast S_3$.

Recall that the modular group $\Gamma = C_2 \ast C_3 = \langle s,t | s^2=1 = t^3 \rangle$ is a quotient of $B_3$ by dividing out the central element $c = (\sigma_1 \sigma_2)^3$, so that we can take $t = \overline{\sigma}_1 \overline{\sigma}_2$ and $s = \overline{\sigma}_1 \overline{\sigma}_2 \overline{\sigma}_1$. Hence, any irreducible $n$-dimensional representation $\phi : B_3 \rTo GL_n$ will be isomorphic to one of the form
\[
\phi(\sigma_1) = \lambda \psi(\overline{\sigma}_1), \quad \text{and} \quad \phi(\sigma_2) = \lambda \psi(\overline{\sigma}_2) \]
for some $\lambda \in \C^*$ and $\psi : \Gamma \rTo GL_n$ an $n$-dimensional irreducible representation of $\Gamma = \langle s,t \rangle = \langle \overline{\sigma}_1,\overline{\sigma_2} \rangle$. With $S_3 = \langle s_1,s_2 | s_1 s_2 s_1 = s_2 s_1 s_2, s_1^2=1=s_2^2 \rangle$, we consider the amalgamated free product
\[
G = \Gamma \ast_{C_3} S_3 \]
in which  the generator of $C_3$ is equal to $t = \overline{\sigma}_1 \overline{\sigma}_2$ in $\Gamma$ and to $s_1 s_2$ in $S_3$.

We will impose conditions on $\psi$ such that it extends to a (necessarily irreducible) representations of $G$. Then, if this is possible, as $\psi(\overline{\sigma}_1 \overline{\sigma}_2) = \psi(s_1 s_2)$ and as the defining equations $(1)$,$(4)$ and $(5)$ of $LB_3$ are homogeneous in the $\sigma_i$ it will follow that 
\[
\phi(\sigma_i) = \lambda \psi(\overline{\sigma}_i), \quad \text{and} \quad \phi(s_i) = \psi(s_i) \]
is a representation of $LB_3$ extending the irreducible representation $\phi$ of $B_3$.

\section{The result}

Bruce Westbury has shown in \cite{Westbury} that the variety $\wis{iss}_n~\Gamma$ classifying isomorphism classes of $n$-dimensional semi-simple $\Gamma$-representations decomposes as a disjoint union of irreducible components
\[
\wis{iss}_n~\Gamma = \bigsqcup_{\alpha} \wis{iss}_{\alpha}~\Gamma \]
where $\alpha=(a,b;x,y,z) \in \mathbb{N}^{\oplus 5}$ satisfying $a+b=n=x+y+z$. Moreover, if $xyz \not= 0$ then the component $\wis{iss}_{\alpha}~\Gamma$ contains a Zariski open and dense subset of irreducible $\Gamma$-representations if and only if $max(x,y,z) \leq min(a,b)$. In this case, the dimension of $\wis{iss}_{\alpha}~\Gamma$ is equal to $1+n^2-(a^2+b^2+x^2+y^2+z^2)$. In going from irreducible $\Gamma$-representations to irreducible $B_3$-representations we multiply by $\lambda \in \C^*$. As a result, it is shown in \cite{Westbury} that there is a $\pmb{\mu}_6$-action on the components $\wis{iss}_{\alpha}~\Gamma$ leading to the same component of $B_3$-representations. That is, the variety $\wis{irr}_n~B_3$ classifying isomorphism classes of irreducible $n$-dimensional $B_3$-representations decomposes into irreducible components
\[
\wis{irr}_n~B_3 = \bigcup_{\alpha} \wis{irr}_{\alpha}~B_3 \]
where $\alpha=(a,b;x,y,z)$ satisfies $a+b=n=x+y+z$, $a \geq b \geq x=max(x,y,z)$.

\begin{theorem} A Zariski open and dense subset of irreducible $\Gamma$-representations in $\wis{iss}_{\alpha}~\Gamma$ extends to the group $G = \Gamma \ast_{C_3} S_3$ if and only if there are natural numbers $u,v,w$ with $w \geq max(u,v)$ such that
\[
\alpha = (v+w,u+w;u+v,w,w) \]
As a consequence, a Zariski open and dense subset of irreducible $B_3$-representations in $\wis{irr}_{\alpha}~B_3$ extends to the three-component loop braid group $LB_3$ if there are natural numbers $u \leq v \leq w$ such that $\alpha=(a,b;x,y,z)$ with $x = max(x,y,z)$ and
\[
a=v+w,~b=u+w,~\{ x,y,z \} = \{ u+v,w,w \} \]
\end{theorem}

Observe that the first dimension $n$ allowing an admissible $5$-tuple not satisfying this condition is $n=6$ with $\alpha=(3,3;3,2,1)$.

\section{The proof}

If $V$ is an $n$-dimensional $G = \Gamma \ast_{C_3} S_3 \simeq C_2 \ast S_3$-representation, then by restricting it to the subgroups $C_2$ and $S_3$ we get decomposition of $V$ into
\[
S_+^{\oplus a} \oplus S_-^{\oplus b} = V \downarrow_{C_2} = V = V \downarrow_{S_3} = T^{\oplus x} \oplus S^{\oplus y} \oplus P^{\oplus z} \]
where $\{ S_+,S_- \}$ are the $1$-dimensional irreducibles of $C_2$, $T$ is the trivial $S_3$-representation, $S$ the sign representation and $P$ the $2$-dimensional irreducible $S_3$-representation. Clearly we must have $a+b = n = x+y + 2z$ and once we choose bases in each of these irreducibles we have that $V$ itself determines a representation of the following quiver setting
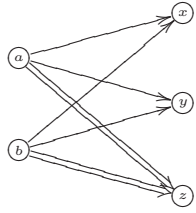
\[
\xymatrix@=.3cm{
& & & & \vtx{x} \\
\vtx{a} \ar[rrrru] \ar[rrrrd] \ar@{=>}[rrrrddd] & & & & \\
& & & & \vtx{y} \\
\vtx{b} \ar[rrrru] \ar[rrrruuu] \ar@{=>}[rrrrd] & & & & \\
& & & &  \vtx{z} }
\]
where the arrows give the block-decomposition of the base-change matrix $B$ from the chosen basis of $V \downarrow_{C_2}$ to the chosen basis of $V \downarrow_{S_3}$. Isomorphism classes of irreducible $G$-representations correspond to isomorphism classes of $\theta$-stable quiver representation of dimension vector $\beta=(a,b;x,y,z)$ for the stability structure $\theta=(-1,-1;1,1,2)$.  The minimal dimension vectors of $\theta$-stable representations are
\[
\begin{cases}
\alpha_1 = (1,0;1,0,0) \\
\alpha_2 = (1,0;0,1,0) \\
\alpha_3 = (0,1;1,0,0) \\
\alpha_4 = (0,1;0,1,0) \\
\alpha_5 = (1,1;0,0,1)
\end{cases}
\]
which give us unique $1$-dimensional $G$-representations $S_1,S_2,S_3,S_4$ and a $2$-parameter family of $2$-dimensional irreducible $G$-representations from which we choose $S_5$. By the results of \cite{Adri-LB}, the local structure of the component $\wis{iss}_{\beta}~G$ for $\beta=(p+q+t,r+s+t;p+r,q+s,t)$ in a neighborhood of the semi-simple $G$-representation
\[
M = S_1^{\oplus p} \oplus S_2^{\oplus q} \oplus S_3^{\oplus r} \oplus S_4^{\oplus s} \oplus S_5^{\oplus t} \]
is \'etale equivalent to the local structure of the quiver-quotientvariety of the setting below at the zero-representation
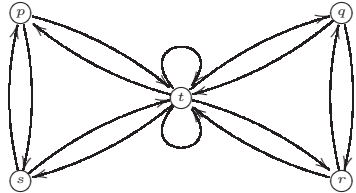
\[
\xymatrix@=.8cm{
\vtx{p} \ar@/^1ex/[rrd] \ar@/^1ex/[dd] & & &  & \vtx{q} \ar@/^1ex/[lld] \ar@/^1ex/[dd] \\
& &  \vtx{t} \ar@/^1ex/[llu] \ar@/^1ex/[rru] \ar@/^1ex/[lld] \ar@/^1ex/[rrd] \ar@(ul,ur) \ar@(dl,dr) & & \\
\vtx{s} \ar@/^1ex/[uu] \ar@/^1ex/[rru] & & & & \vtx{r} \ar@/^1ex/[uu] \ar@/^1ex/[llu] }
\]
Hence, $\wis{iss}_{\beta}~G$ will contain a Zariski open and dense subset of irreducible representations if and only if $\gamma = (p,q,r,s,t)$ is a simple dimension vector for this quiver, which by \cite{LBProcesi} is equivalent to $\gamma$ being either $(1,0,0,1,0)$ or $(0,1,1,0,0)$ or satisfying the inequalities
\[
p \leq s+t,~q \leq r+t,~r \leq q+t,~s \leq p+t \]
Having determined the components containing irreducible $G$-representations, we have to determine those containing a Zariski open subset which remain irreducible when restricted to $\Gamma$. 

As $\Gamma = C_2 \ast C_3$ any $\Gamma$-representation $V$ corresponds to a semi-stable quiver represenation for the setting
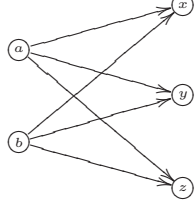
\[
\xymatrix@=.3cm{
& & & & \vtx{x} \\
\vtx{a} \ar[rrrru] \ar[rrrrd] \ar[rrrrddd] & & & & \\
& & & & \vtx{y} \\
\vtx{b} \ar[rrrru] \ar[rrrruuu] \ar[rrrrd] & & & & \\
& & & &  \vtx{z} }
\]
when
\[
V \downarrow_{C_2} = S_+^{\oplus a} \oplus S_-^{\oplus b} \quad \text{and} \quad V \downarrow_{C_3} = T_1^{\oplus x} \oplus T_{\rho}^{\oplus y} \oplus T_{\rho^2}^{\oplus z} \]
with $\{ T_1,T_{\rho}, T_{\rho^2} \}$ the irreducible $C_3$-representations. Because $T \downarrow_{C_3} = T_1 = S \downarrow_{C_3}$ and $P \downarrow_{C_3} = T_{\rho} \oplus T_{\rho^2}$ we have that $M \downarrow_{\Gamma}$ has dimension vector
\[
\alpha = (a,b;x,y,z) = (p+q+t,r+s+t;p+q+r+s,t,t) \]
which satisfies the condition that $max(x,y,z) \leq min(a,b)$ if and only if $t \geq r+s$ and $t \geq p+q$. Setting $u=r+s$, $v=p+q$ and $w=t$, the statement of Theorem~1 follows.


\begin{thebibliography}{10}

\bibitem{Adri-LB}
Jan Adriaenssens and Lieven Le Bruyn, {\it Local quivers and stable representations}, Communications in Algebra, 31 (2003) 1777-1797


\bibitem{Baez}
John Baez, Derek Wise and Alissa Crans, {\it Exotic statistics for strings in $4D$ BF theory}, Adv. Theor. Math. Phys. {\bf 11} (2007) 707-749, {\tt arXiv:gr-qc/0603085}

\bibitem{Bruillard}
Paul Bruillard, Seung-Moon Hong, Julia Yael Plavnik, Eric C. Rowell, Liang Chang and Michael Yuan Sun, {\it Low-dimensional representations of the three component loop braid group}, {\tt arXiv:1508.00005} (2015)


\bibitem{LBbraid1}
Lieven Le Bruyn, {\it Dense families of $B_3$-representations and braid reversion}, Journal of Pure and Appl. Algebra 215 (2011) 1003-1014, {\tt math.RA/1003.1610}


\bibitem{LBProcesi}
Lieven Le Bruyn and Claudio Procesi, {\it Semisimple representations of quivers}, Trans. Amer. Math. Soc. {\bf 317} (1990) 585-598 


\bibitem{TubaWenzl}
Imre Tuba and Hans Wenzl, {\it Representations of the braid group $B_3$ and of $SL(2,\Z)$}, Pacific J. Math. 197 (2001)

\bibitem{Westbury}
Bruce Westbury, {\it On the character varieties of the modular group}, preprint Nottingham University (1995)

\end{thebibliography}
\end{document}